\documentclass[ejs,preprint]{imsart}

\pdfoutput=1
\usepackage{amsfonts,amssymb,amsmath,amscd,amsthm,latexsym,color}
\usepackage{natbib}
\bibpunct{(}{)}{,}{a}{,}{,}
\usepackage{amsfonts,amssymb,amsmath,amscd,amsthm,latexsym}
\usepackage{natbib}
\bibpunct{(}{)}{,}{a}{,}{,}

\usepackage{graphicx}
\usepackage{fullpage}

\newtheorem{example}{Example}

\newsavebox{\fmbox}

\renewcommand{\vss}{\vspace{0.25cm}}

\begin{document}

\begin{frontmatter}

\title{On resolving the Savage--Dickey paradox}
\runtitle{On the Savage--Dickey paradox}

\begin{aug}
\author{Jean-Michel Marin\ead[label=e1]{jean-michel.marin@univ-montp2.fr}}
\and
\author{Christian P.~Robert\ead[label=e2]{xian@ceremade.dauphine.fr}}

\runauthor{Robert, C.P. and Marin, J.-M.}

\affiliation{Institut de Math\'ematiques et Mod\'elisation de Montpellier,\\ 
Universit\'e Montpellier 2, and CREST, INSEE, Paris}
\affiliation{Universit\'e Paris-Dauphine, CEREMADE, and CREST, Paris}
\address{Institut de Math\'ematiques et Mod\'elisation de Montpellier,\\
Universit\'e Montpellier 2, Case Courrier 51\\ 34095 Montpellier cedex 5, France,\\ \printead{e1}}
\address{Universit\'e Paris-Dauphine, CEREMADE\\ 75775 Paris cedex 16, France,\\
         CREST\\ 92245 Malakoff cedex, France\\ \printead{e2}}

\end{aug}

\begin{abstract}
When testing a null hypothesis $H_0:~\theta=\theta_0$ in a Bayesian framework,
the Savage--Dickey ratio \citep{dickey:1971} is known as a specific representation 
of the Bayes factor \citep{ohagan:forster:2004} that only uses the posterior distribution 
under the alternative hypothesis at $\theta_0$, thus allowing for a
plug-in version of this quantity. 
We demonstrate here that the Savage--Dickey representation is in fact a generic representation
of the Bayes factor and that it fundamentally relies on specific measure-theoretic versions of the densities 
involved in the ratio, instead of being a special identity imposing some mathematically void constraints 
on the prior distributions.
We completely clarify the measure-theoretic foundations of the Savage--Dickey representation as well as 
of the later generalisation of \cite{verdinelli:wasserman:1995}.  We provide furthermore
a general framework that produces a converging approximation of the Bayes factor that is unrelated
with the approach of \cite{verdinelli:wasserman:1995} and propose a comparison of this new
approximation with their version, as well as with bridge sampling and Chib's approaches.
\end{abstract}


\end{frontmatter}

\section{Introduction}

From a methodological viewpoint, testing a null hypothesis $H_0:~x\sim f_0(x|\omega_0)$ versus the alternative
$H_a:~x\sim f_1(x|\omega_1)$  in a Bayesian framework requires the introduction of two prior distributions, 
$\pi_0(\omega_0)$ and $\pi_1(\omega_1)$, that are defined on the respective parameter spaces.
In functional terms, the core object of the Bayesian approach to testing and model choice, the Bayes factor \citep{jeffreys:1939,
robert:2001,ohagan:forster:2004}, is indeed a ratio of two marginal densities taken at the same observation $x$,
$$
B_{01}(x) = \dfrac{\int \pi_0(\omega_0) f_0(x|\omega_0)\,\text{d}\omega_0}
{\int \pi_1(\omega_1) f_1(x|\omega_1)\,\text{d}\omega_1} = \dfrac{m_0(x)}{m_1(x)}\,.
$$
(This quantity $B_{01}(x)$ is then compared to $1$ in order to decide about the strength of the support of the data in
favour of $H_0$ or $H_a$.)
It is thus mathematically clearly and uniquely defined, provided both integrals exist and differ from both $0$ and
$\infty$. The practical computation of the Bayes factor has generated a large literature on approximative \citep[see, 
e.g.][]{chib:1995,gelman:meng:1998,chen:shao:ibrahim:2000,chopin:robert:2010}, seeking improvements in numerical precision.

The Savage--Dickey \citep{dickey:1971} representation of the Bayes factor
is primarily known as a special identity that relates the Bayes factor to the 
posterior distribution which corresponds to the more complex hypothesis. As described in \cite{verdinelli:wasserman:1995} and
\citeauthor{chen:shao:ibrahim:2000} (2000, pages 164-165), this representation has practical implications as a basis for
simulation methods. However, as stressed in \cite{dickey:1971} and \cite{ohagan:forster:2004}, 
the foundation of the Savage--Dickey representation is clearly theoretical. 

More specifically, when considering a testing problem with
an embedded model, $H_0:\theta=\theta_0$, and a nuisance parameter $\psi$, i.e.~when $\omega_1$ can be decomposed
as $\omega_1=(\theta,\psi)$ and when $\omega_0=(\theta_0,\psi)$, for a sampling distribution $f(x|\theta,
\psi)$, the plug-in representation 
\begin{equation}\label{eq:dickey}
B_{01}(x) = \dfrac{\pi_1(\theta_0|x)}{\pi_1(\theta_0)}\,,
\end{equation}
with the obvious notations for the marginal distributions
$$
\pi_1(\theta) = \int \pi_1(\theta,\psi)\text{d}\psi 
\quad\text{and}\quad 
\pi_1(\theta|x) = \int \pi_1(\theta,\psi|x)\text{d}\psi\,,
$$
holds under Dickey's (1971) assumption that the conditional prior density of $\psi$ 
under the alternative model, given $\theta=\theta_0$, $\pi_1(\psi|\theta_0)$, 
is equal to the prior density under the null hypothesis, $\pi_0(\psi)$,
\begin{equation}\label{eq:savage}
\pi_1(\psi|\theta_0)=\pi_0(\psi)\,.
\end{equation}
Therefore, Dickey's (1971) identity \eqref{eq:dickey} reduces the Bayes factor to the ratio of
the posterior over the prior marginal densities of $\theta$ under the alternative model, 
taken at the tested value $\theta_0$. The Bayes factor is thus expressed as an amount of information
brought by the data and this helps in its justification as a model choice tool. (See also 
\citealp{consonni:veronese:2008}.)

In order to illustrate the Savage--Dickey representation, consider the artificial example of computing
the Bayes factor between the models
$$
\mathfrak{M}_0:\quad x|\psi\sim\mathcal{N}(\psi,1),\quad \psi\sim\mathcal{N}(0,1)\,,
$$
and
$$
\mathfrak{M}_1:\quad x|\theta,\psi\sim\mathcal{N}(\psi,\theta),\quad \psi|\theta\sim\mathcal{N}(0,\theta),\quad
\theta\sim I\mathcal{G}(1,1)\,,
$$
which is equivalent to testing the null hypothesis $H_0:\theta=\theta_0=1$ against the alternative
$H_1:\theta\neq 1$ when $x|\theta,\psi\sim\mathcal{N}(\psi,\theta)$. In that case,
model $\mathfrak{M}_0$ clearly is embedded in model $\mathfrak{M}_1$.
We have
$$
m_0(x)=\exp\left(-x^2/4\right)\big/ (\sqrt{2}\sqrt{2\pi})\quad\mbox{and}\quad
m_1(x)=\left(1+x^2/4\right)^{-3/2}\Gamma(3/2)\big/ (\sqrt{2}\sqrt{2\pi})\,,
$$
and therefore
$$
B_{01}(x)=\Gamma(3/2)^{-1}\left(1+x^2/4\right)^{3/2}\exp\left(-x^2/4\right)\,.
$$
Dickey's assumption \eqref{eq:savage} on the prior densities is satisfied, since
$$
\pi_1(\psi|\theta_0)=\frac{1}{\sqrt{2\pi}}\exp\left(-\psi^2/2\right)=\pi_0(\psi)\,.
$$
Therefore, since
$$
\pi_1(\theta)=\theta^{-2}\exp\left(-\theta^{-1}\right)\,,\quad\pi_1(\theta_0)=\exp(-1)\,,
$$
and
\begin{align*}
\pi_1(\theta|x)&=\Gamma(3/2)^{-1}\left(1+x^2/4\right)^{3/2}\theta^{-5/2}
\exp\left(-\theta^{-1}\left(1+x^2/4\right)\right)\mathbb{I}_{\theta>0}\,,\\
\pi_1(\theta_0|x)&=\Gamma(3/2)^{-1}\left(1+x^2/4\right)^{3/2}\exp\left(-\left(1+x^2/4\right)\right)\,,
\end{align*}
we clearly recover the Savage--Dickey representation
$$
B_{01}(x)=\Gamma(3/2)^{-1}\left(1+x^2/4\right)^{3/2}\exp\left(-x^2/4\right)=\pi_1(\theta_0|x)/\pi_1(\theta_0)\,.
$$

While the difficulty with the representation \eqref{eq:dickey} is usually addressed in terms of computational aspects, 
given that $\pi_1(\theta|x)$ is rarely available in closed form, we argue in the current paper
that the Savage--Dickey representation faces challenges of a deeper nature that led us to consider it a `paradox'. 
First, 
by considering both prior and posterior marginal distributions of $\theta$ uniquely {\em under the alternative model,} 
\eqref{eq:dickey} seems to indicate that the posterior probability of the null hypothesis $H_0:\theta=\theta_0$ 
is contained within the alternative hypothesis posterior distribution, even though the set of $(\theta,\psi)$'s 
such that $\theta=\theta_0$ has a zero probability under this alternative distribution. 
Second, 
as explained in Section \ref{sec:measure}, an even more fundamental difficulty with assumption \eqref{eq:savage} 
is that it is meaningless when examined (as it should) within the mathematical axioms of measure theory.

Having stated those mathematical difficulties with the Savage--Dickey representation, we proceed 
to show in Section \ref{sec:montecarl} that similar identities hold
under no constraint on the prior distributions. In Section \ref{sec:montecarl},
we derive computational algorithms that exploit these representations to approximate the Bayes factor, in
an approach that differs from the earlier solution of \cite{verdinelli:wasserman:1995}.
The paper concludes with an illustration in the setting of variable selection within a probit model.

\section{A measure-theoretic paradox}\label{sec:measure}

When considering a standard probabilistic setting where the dominating measure on the parameter space is the Lebesgue measure, 
rather than a counting measure, the conditional
density $\pi_1(\psi|\theta)$ is rigorously \citep{billingsley:1986} defined as the density of the conditional probability distribution
or, equivalently, by the condition that
$$
\mathbb{P}((\theta,\psi)\in A_1\times A_2) = \int_{A_1} \int_{A_2} \pi_1(\psi|\theta) \,\text{d}\psi\,\pi_1(\theta)\,\text{d}\theta
		                           = \int_{A_1\times A_2} \pi_1(\theta,\psi) \text{d}\psi\,\text{d}\theta\,,
$$
for all measurable sets $A_1\times A_2$, when $\pi_1(\theta)$ is the associated marginal density of $\theta$. 
Therefore, this identity points out the well-known fact that the conditional density function $\pi_1(\psi|\theta)$ is
defined up to a set of measure zero both in $\psi$ for {\em every} value of $\theta$ {\em and} in $\theta$. 
This implies that changing arbitrarily the value of the
{\em function} $\pi_1(\cdot|\theta)$ for a negligible collection of values of $\theta$ does not impact the properties of the
conditional distribution. 

In the setting where the Savage--Dickey representation is advocated, the value $\theta_0$ to be tested is not determined from
the observations but it is instead given in advance since this is a testing problem. Therefore the density function
$$
\pi_1(\psi|\theta_0)
$$
may be chosen in a {\em completely arbitrary} manner and there is no possible reason for a unique representation 
of $\pi_1(\psi|\theta_0)$ that can be found within measure theory. This implies that there always
is a version of the conditional density $\pi_1(\psi|\theta_0)$ such that  Dickey's (1971) condition \eqref{eq:savage} is 
satisfied---as well as, conversely, there are an infinity of versions for which it is {\em not} satisfied---. As a result, 
from a mathematical perspective, condition \eqref{eq:savage} cannot be seen as an {\em assumption} on the prior $\pi_1$ 
without further conditions, contrary to what is stated in the original \cite{dickey:1971} and later in \cite{ohagan:forster:2004}, \cite{consonni:veronese:2008}
and \cite{wetzels:grasman:wagenmakers:2010}.
This difficulty is the  first part of what we call the {\em Savage--Dickey paradox}, namely that, as stated, the representation
\eqref{eq:dickey} relies on a mathematically void constraint on the prior distribution. In the specific case of the artificial
example introduced above, the choice of the conditional density $\pi_1(\psi|\theta_0)$ is therefore arbitrary: if we pick for
this density the density of the $\mathcal{N}(0,1)$ distribution, there is agreement between $\pi_1(\psi|\theta_0)$ and 
$\pi_0(\psi)$, while, if we select instead the function $\exp(+\psi^2/2)$, which is not a density, there is no agreement in 
the sense of condition \eqref{eq:savage}. The paradox is that this disagreement has no consequence whatsoever in the
Savage--Dickey representation.

The second part of the Savage--Dickey paradox is that the representation \eqref{eq:dickey} is solely valid for a specific and unique 
choice of a version of the density for both the conditional density $\pi_1(\psi|\theta_0)$ and the joint density $\pi_1(\theta_0,\psi)$. 
When looking at the derivation of \eqref{eq:dickey}, the choices of some specific versions 
of those densities are indeed noteworthy: in the following development, 
\begin{alignat*}{2}
B_{01}(x) & = \dfrac{\int \pi_0(\psi) f(x|\theta_0,\psi)\,\text{d}\psi}{\int \pi_1(\theta,\psi)f(x|\theta,\psi)\,
              \text{d}\psi\text{d}\theta} & \qquad & \text{[by definition]} \\ 
          & = \dfrac{\int \pi_1(\psi|\theta_0) f(x|\theta_0,\psi)\,\text{d}\psi\,\pi_1(\theta_0)}
              {\int \pi_1(\theta,\psi)f(x|\theta,\psi)\,\text{d}\psi\text{d}\theta\,\pi_1(\theta_0)}
          && \text{[using a specific version of $\pi_1(\psi|\theta_0)$]}\\
          & = \dfrac{\int \pi_1(\theta_0,\psi) f(x|\theta_0,\psi)\,\text{d}\psi}{m_1(x)\pi_1(\theta_0)}
          && \text{[using a specific version of $\pi_1(\theta_0,\psi)$]}\\
          & = \dfrac{\pi_1(\theta_0|x)}{\pi_1(\theta_0)}\,,&&\text{[using a specific version of $\pi_1(\theta_0|x)$]}
\end{alignat*}
the second equality depends on a specific choice of the version of $\pi_1(\psi|\theta_0)$ 
but not on the choice of the version of $\pi_1(\theta_0)$,
while the third equality depends on a specific choice of the version of $\pi_1(\psi,\theta_0)$ as
equal to $\pi_0(\psi) \pi_1(\theta_0)$, thus related to the choice of the version of $\pi_1(\theta_0)$. 
The last equality leading to the Savage--Dickey representation relies on 
the choice of a specific version of $\pi_1(\theta_0|x)$ as well, namely that the constraint
$$
\dfrac{\pi_1(\theta_0|x)}{\pi_1(\theta_0)} = 
\dfrac{\int \pi_0(\psi) f(x|\theta_0,\psi)\,\text{d}\psi}{m_1(x)}
$$
holds, where the right hand side is equal to the Bayes factor $B_{01}(x)$ and is therefore independent from the
version. This rigorous analysis implies that the Savage--Dickey representation is tautological, due to the 
availability of a version of the posterior density that makes it hold.

As an illustration, consider once again the artificial example above. As already stressed, the value to be tested $\theta_0=1$ is 
set prior to the experiment.  Thus, without modifying either the prior distribution under model $\mathfrak{M}_1$ 
or the marginal posterior distribution of the parameter $\theta$ under model $\mathfrak{M}_1$, and in a completely 
rigorous measure-theoretic framework, we can select
$$
\pi_1(\theta_0)=100=\pi_1(\theta_0|x)\,.
$$  
For that choice, we obtain
$$
\pi_1(\theta_0|x)/\pi_1(\theta_0) = 1 \neq B_{01}(x)= \Gamma(3/2)^{-1}\left(1+x^2/4\right)^{3/2}\exp\left(-x^2/4\right)\,.
$$
Hence, for this specific choice of the densities, the Savage--Dickey representation does not hold.

\vspace{0.5cm} \cite{verdinelli:wasserman:1995} have proposed a generalisation of the Savage--Dickey density ratio when
the constraint (\ref{eq:savage}) on the prior densities is not verified (we stress again that this is a mathematically
void constraint on the respective prior distributions). \cite{verdinelli:wasserman:1995} state that
\begin{alignat*}{2}
B_{01}(x) & = \dfrac{\int \pi_0(\psi) f(x|\theta_0,\psi)\,\text{d}\psi}{m_1(x)} & \qquad & \text{[by definition]}\\ 
          & = \pi_1(\theta_0|x)\dfrac{\int \pi_0(\psi) f(x|\theta_0,\psi)\,\text{d}\psi}{m_1(x)\pi_1(\theta_0|x)} 
		& \qquad & \text{[for any version of $\pi_1(\theta_0|x)$]}\\
          & = \pi_1(\theta_0|x)\int\dfrac{\pi_0(\psi) f(x|\theta_0,\psi)}{m_1(x)\pi_1(\theta_0|x)}
		\dfrac{\pi_1(\psi|\theta_0)}{\pi_1(\psi|\theta_0)}\,\text{d}\psi & \qquad & \text{[for any version of $\pi_1(\psi|\theta_0)$]}\\
          & = \pi_1(\theta_0|x)\int\dfrac{\pi_0(\psi)}{\pi_1(\psi|\theta_0)}\,\dfrac{f(x|\theta_0,\psi) 
		\pi_1(\psi|\theta_0)\,\text{d}\psi}{m_1(x)\pi_1(\theta_0|x)}\,\dfrac{\pi_1(\theta_0)}{\pi_1(\theta_0)} 
		& \qquad & \text{[for any version of $\pi_1(\theta_0)$]}\\
          & = \dfrac{\pi_1(\theta_0|x)}{\pi_1(\theta_0)}\,\int\dfrac{\pi_0(\psi)}{\pi_1(\psi|\theta_0)}\, \pi_1(\psi|\theta_0,x)\,\text{d}\psi
		& \qquad & \text{[for a specific version of $\pi_1(\psi|\theta_0,x)$]}\\
          & = \dfrac{\pi_1(\theta_0|x)}{\pi_1(\theta_0)}\,\mathbb{E}^{\pi_1(\psi|x,\theta_0)}
	      \left[\dfrac{\pi_0(\psi)}{\pi_1(\psi|\theta_0)}\right]\,.
\end{alignat*}
This representation of \cite{verdinelli:wasserman:1995} therefore remains valid for any choice of versions 
for $\pi_1(\theta_0|x)$, $\pi_1(\theta_0)$, $\pi_1(\psi|\theta_0)$, provided the conditional density $\pi_1(\psi|\theta_0,x)$ is
defined by
$$
\pi_1(\psi|\theta_0,x) = \dfrac{f(x|\theta_0,\psi) \pi_1(\psi|\theta_0)\pi_1(\theta_0)}{m_1(x)\pi_1(\theta_0|x)}\,,
$$
which obviously means that the Verdinelli--Wasserman representation 
\begin{equation}
B_{01}(x) = \dfrac{\pi_1(\theta_0|x)}{\pi_1(\theta_0)}\,\mathbb{E}^{\pi_1(\psi|x,\theta_0)}
\left[\dfrac{\pi_0(\psi)}{\pi_1(\psi|\theta_0)}\right]
\end{equation}
is dependent on the choice of a version of $\pi_1(\theta_0)$.

\vspace{0.5cm} We now establish that an alternative representation of the Bayes factor is available and can be exploited towards
approximation purposes. When considering the Bayes factor
$$
B_{01}(x) = \dfrac{\int \pi_0(\psi) f(x|\theta_0,\psi)\,\text{d}\psi}
{\int \pi_1(\theta,\psi)f(x|\theta,\psi)\, \text{d}\psi\text{d}\theta}\,
\dfrac{\pi_1(\theta_0)}{\pi_1(\theta_0)}\,,
$$
where the right hand side obviously is independent of the choice of the version of $\pi_1(\theta_0)$, 
the numerator can be seen as involving a specific version in $\theta=\theta_0$ of the marginal posterior density
$$
\tilde\pi_1(\theta|x) \propto \int \pi_0(\psi) f(x|\theta,\psi) \,\text{d}\psi\,\pi_1(\theta)\,,
$$
which is associated with the alternative prior $\tilde\pi_1(\theta,\psi)=\pi_1(\theta)\pi_0(\psi)$.
Indeed, this density $\tilde\pi_1(\theta|x)$ appears as the marginal posterior density of the posterior 
distribution defined by the density
$$
\tilde\pi_1(\theta,\psi|x) = \dfrac{ \pi_0(\psi) \pi_1(\theta) f(x|\theta,\psi) }{\tilde{m}_1(x)}\,,
$$
where $\tilde m_1(x)$ is the proper normalising constant of the joint posterior density.
In order to guarantee a Savage--Dickey-like representation of the Bayes factor,
the appropriate version of the marginal posterior density in $\theta=\theta_0$,
$\tilde\pi_1(\theta_0|x)$,  is obtained by imposing
\begin{equation}\label{eq:psudopost}
\dfrac{\tilde\pi_1(\theta_0|x)}{\pi_0(\theta_0)} = \dfrac{\int \pi_0(\psi) 
f(x|\theta_0,\psi) \,\text{d}\psi}{\tilde{m}_1(x)}\,,
\end{equation}
where, once again, the right hand side of the equation is uniquely defined. This constraint amounts to imposing that Bayes' 
theorem holds in $\theta=\theta_0$ instead of almost everywhere (and thus not necessarily in $\theta=\theta_0$).
It then leads to the alternative representation
$$
B_{01}(x) = \dfrac{\tilde\pi_1(\theta_0|x)}{\pi_1(\theta_0)}\,\dfrac{\tilde{m}_1(x)}{m_1(x)}\,,
$$
which holds for any value chosen for $\pi_1(\theta_0)$ provided condition \eqref{eq:psudopost} applies.

This new representation may seem to be only formal, since both $m_1(x)$ and $\tilde m_1(x)$
are usually unavailable in closed form, but we can take advantage of the fact that the 
bridge sampling identity of \cite{torrie:valleau:1977} (see also 
\citealp{gelman:meng:1998}) gives an unbiased estimator of $\tilde m_1(x)/{m}_1(x)$ since
$$
\mathbb{E}^{\pi_1(\theta,\psi|x)} \left[ \dfrac{\pi_0(\psi)\pi_1(\theta) f(x|\theta,\psi) }
{\pi_1(\theta,\psi)f(x|\theta,\psi) } \right] =  
\mathbb{E}^{\pi_1(\theta,\psi|x)} \left[ \dfrac{\pi_0(\psi)}{\pi_1(\psi|\theta)}
\right] = \dfrac{\tilde m_1(x)}{{m}_1(x)}\,.
$$
In conclusion, we obtain the representation
\begin{equation}
B_{01}(x) = \dfrac{\tilde\pi_1(\theta_0|x)}{\pi_1(\theta_0)}\,
\mathbb{E}^{\pi_1(\theta,\psi|x)} \left[ \dfrac{\pi_0(\psi)}{\pi_1(\psi|\theta)}\right]\,,
\label{eq:mr09}
\end{equation}
whose expectation part is uniquely defined (in that it does not depend on the choice of a version 
of the densities involved therein), while the first ratio must satisfy condition \eqref{eq:psudopost}.
We further note that this representation clearly differs from Verdinelli and Wasserman's 
(\citeyear{verdinelli:wasserman:1995}) representation:
\begin{equation}
B_{01}(x)=\dfrac{\pi_1(\theta_0|x)}{\pi_1(\theta_0)}\,\mathbb{E}^{\pi_1(\psi|x,\theta_0)}\left[\dfrac{\pi_0(\psi)}{\pi_1(\psi|\theta_0)}\right]\,,
\label{eq:vw05}
\end{equation}
since \eqref{eq:vw05} uses a specific version of the marginal posterior density on $\theta$ in $\theta_0$, as
well as a specific version of the full conditional posterior density of $\psi$ given $\theta_0$

\section{Computational solutions}\label{sec:montecarl}
In this Section, we consider the computational implications of the above
representation in the specific case of latent variable models, namely under the
practical possibility of a data completion by a latent variable $z$ such that 
$$
f(x|\theta,\psi) = \int f(x|\theta,\psi,z)f(z|\theta,\psi)\,\text{d}z
$$
when $\pi_1(\theta|x,\psi,z) \propto \pi_1(\theta) f(x|\theta,\psi,z)$ is 
available in closed form, including the normalising constant.

\vspace{0.5cm} We first consider a computational solution that approximates the Bayes factor based on our 
novel representation (\ref{eq:mr09}). Given a sample $(\bar \theta^{(1)},\bar \psi^{(1)},\bar z^{(1)}),\allowbreak \ldots,
\allowbreak(\bar \theta^{(T)},\bar \psi^{(T)},\bar z^{(T)})$ 
simulated from (or converging to) the augmented posterior distribution $\tilde\pi_1(\theta,\psi,z|x)$, the sequence
$$
\dfrac{1}{T}\,\sum_{t=1}^T \tilde\pi_1(\theta_0|x,\bar z^{(t)},\bar \psi^{(t)})
$$
converges to $\tilde\pi_1(\theta_0|x)$ in $T$ under the following constraint on the 
selected version of $\tilde\pi_1(\theta_0|x,z,\psi)$ used therein:
$$
\dfrac{\tilde\pi_1(\theta_0|x,z,\psi)}{\pi_1(\theta_0)}=
\dfrac{f(x,z|\theta_0,\psi)}{\int f(x,z|\theta,\psi) \pi_1(\theta)\,\text{d}\theta}\,.
$$
which again amounts to imposing that Bayes' theorem holds in $\theta=\theta_0$ for $\tilde\pi_1(\theta|x,z,\psi)$
rather than almost everywhere. (Note once more that the right hand side is uniquely defined,
i.e.~that it does not depend on a specific version.)
Therefore, provided iid or MCMC simulations from the joint target $\tilde\pi_1(\theta,\psi,z|x)$ are available, the
converging approximation to the Bayes factor $B_{01}(x)$ is then
$$
\dfrac{1}{T}\sum_{t=1}^T \dfrac{\tilde\pi_1(\theta_0|x,\bar z^{(t)},\bar \psi^{(t)})}{\pi_1(\theta_0)}\,
\dfrac{\tilde{m}_1(x)}{m_1(x)}\,.
$$
(We stress that the simulated sample is produced for the artificial target $\tilde\pi_1(\theta,\psi,z|x)$ rather than
the true posterior $\pi_1(\theta,\psi,z|x)$ if $\tilde\pi_1(\theta,\psi)\ne\pi_1(\theta,\psi)$.)
Moreover, if $(\theta^{(1)},\psi^{(1)}), \allowbreak \ldots, \allowbreak(\theta^{(T)},\psi^{(T)})$ is a sample independently
simulated from (or converging to) $\pi_1(\theta,\psi|x)$, then
$$
\dfrac{1}{T}\,\sum_{t=1}^T \dfrac{\pi_0(\psi^{(t)})}{\pi_1(\psi^{(t)}|\theta^{(t)})}
$$
is a convergent and unbiased estimator of $\tilde{m}_1(x)/m_1(x)$.
Therefore, the computational solution associated to our representation \eqref{eq:mr09} of $B_{01}(x)$
leads to the following unbiased estimator of the Bayes factor:  
\begin{equation}\label{eq:arrox-mr09}
\widehat{B_{01}}^{\text{MR}}(x) = \dfrac{1}{T}\,
\sum_{t=1}^T \dfrac{\tilde\pi_1(\theta_0|x,\bar z^{(t)},\bar\psi^{(t)})}{\pi_1(\theta_0)}\,
\dfrac{1}{T}\,\sum_{t=1}^T \dfrac{\pi_0(\psi^{(t)})}{\pi_1(\psi^{(t)}|\theta^{(t)})}\,.
\end{equation}
Note that
$$
\mathbb{E}^{\tilde\pi_1(\theta,\psi|x)} \left[ \dfrac{\pi_1(\theta,\psi)f(x|\theta,\psi) }
{\pi_0(\psi)\pi_1(\theta) f(x|\theta,\psi) } \right] =  
\mathbb{E}^{\tilde\pi_1(\theta,\psi|x)} \left[ \dfrac{\pi_1(\psi|\theta)}
{\pi_0(\psi)} \right] = \dfrac{m_1(x)}{\tilde{m}_1(x)}
$$
implies that 
$$
T \bigg/ \sum_{t=1}^T \dfrac{\pi_1(\bar \psi^{(t)}|\theta^{(t)})}{\pi_0(\bar \psi^{(t)})}
$$
is another convergent (if biased) estimator of $\tilde{m}_1(x)/m_1(x)$.
The availability of two estimates of the ratio $\tilde{m}_1(x)/m_1(x)$ is a major bonus from a
computational point of view since the comparison of both estimators may allow for the detection of
infinite variance estimators, as well as for coherence of the approximations. The first approach requires
two simulation sequences, one from $\tilde\pi_1(\theta,\psi|x)$ and one from $\pi_1(\theta,\psi|x)$, but this is
a void constraint in that, if $H_0$ is rejected, a sample from the alternative hypothesis posterior
will be required no matter what. Although we do not pursue this possibility in the current paper, note that
a comparison of the different representations (including Verdinelli and Wasserman's, 1995, as exposed below) could
be conducted by expressing them in the bridge sampling formalism \citep{gelman:meng:1998}.

\vspace{0.5cm} We now consider a computational solution that approximates the Bayes factor and is based
on \cite{verdinelli:wasserman:1995}'s representation (\ref{eq:vw05}). Given a sample $(\theta^{(1)},
\psi^{(1)},z^{(1)}),\allowbreak \ldots, \allowbreak(\theta^{(T)},\psi^{(T)},z^{(T)})$ simulated from 
(or converging to) $\pi_1(\theta,\psi,z|x)$, the sequence
$$
\dfrac{1}{T}\,\sum_{t=1}^T \pi_1(\theta_0|x,z^{(t)},\psi^{(t)})
$$
converges to $\pi_1(\theta_0|x)$ under the following constraint on the 
selected version of $\pi_1(\theta_0|x,z,\psi)$ used there:
$$
\dfrac{\pi_1(\theta_0|x,z,\psi)}{\pi_1(\theta_0)}=
\dfrac{f(x,z|\theta_0,\psi)}{\int f(x,z|\theta,\psi) \pi_1(\theta)\,\text{d}\theta}\,.
$$
Moreover, if $\left(\tilde \psi^{(1)},\tilde z^{(1)}\right),\ldots,\left(\tilde \psi^{(T)},\tilde z^{(T)}\right)$ 
is a sample generated from (or converging to) $\pi_1(\psi,z|x,\theta_0)$, the sequence
$$
\frac{1}{T}\,\sum_{t=1}^T\frac{\pi_0(\tilde\psi^{(t)})}{\pi_1(\tilde\psi^{(t)}|\theta_0)}
$$
is converging to 
$$
\mathbb{E}^{\pi_1(\psi|x,\theta_0)}\left[\dfrac{\pi_0(\psi)}{\pi_1(\psi|\theta_0)}\right]
$$
under the constraint 
$$
\pi_1(\psi,z|\theta_0,x) \propto f(x,z|\theta_0,\psi) \pi_1(\psi|\theta_0)\,.
$$
Therefore, the computational solution associated to the \cite{verdinelli:wasserman:1995}'s representation of $B_{01}(x)$
(\ref{eq:vw05}) leads to the following unbiased estimator of the Bayes factor:  
\begin{equation}\label{eq:arrox-vw05}
\widehat{B_{01}}^{\text{VW}}(x) = \dfrac{1}{T}\,\sum_{t=1}^T \dfrac{\pi_1(\theta_0|x,z^{(t)},\psi^{(t)})}{\pi_1(\theta_0)}\,
\dfrac{1}{T}\,\sum_{t=1}^T \dfrac{\pi_0(\tilde\psi^{(t)})}{\pi_1(\tilde\psi^{(t)}|\theta_0)}\,.
\end{equation}
Although, at first sight, the approximations \eqref{eq:arrox-mr09} and \eqref{eq:arrox-vw05}
may look very similar, the simulated sequences used in both approximations differ: the first average
involves simulations from $\tilde\pi_1(\theta,\psi,z|x)$ and from $\pi_1(\theta,\psi,z|x)$, respectively,
while the second average relies on simulations from $\pi_1(\theta,\psi,z|x)$ and from $\pi_1(\psi,z|x,\theta_0)$,
respectively.

\section{An illustration}

Although our purpose in this note is far from advancing the superiority of the Savage--Dickey type representations
for Bayes factor approximation, given the wealth of available solutions for embedded models \citep{chen:shao:ibrahim:2000,
marin:robert:2010}, we briefly consider an example where both Verdinelli and Wasserman's (1995) and our proposal apply. 
The model is the Bayesian posterior distribution of the regression coefficients of a probit model, following the
prior modelling adopted in \cite{marin:robert:2007} that extends \citeauthor{zellner:1986}'s (1971)
$g$-prior to generalised linear models. We take as data the Pima Indian diabetes study available in
R \citep{rmanual} dataset with 332 women registered and build a probit model predicting the presence of diabetes
from three predictors, the glucose concentration, the diastolic blood pressure and the diabetes pedigree function,
assessing the impact of the diabetes pedigree function, i.e.~testing the nullity of the coefficient $\theta$ associated
to this variable. For more details on the statistical and computational issues, see \cite{marin:robert:2010}
since this paper relies on the Pima Indian probit model as benchmark.

\begin{figure}
\includegraphics[width=.6\textwidth]{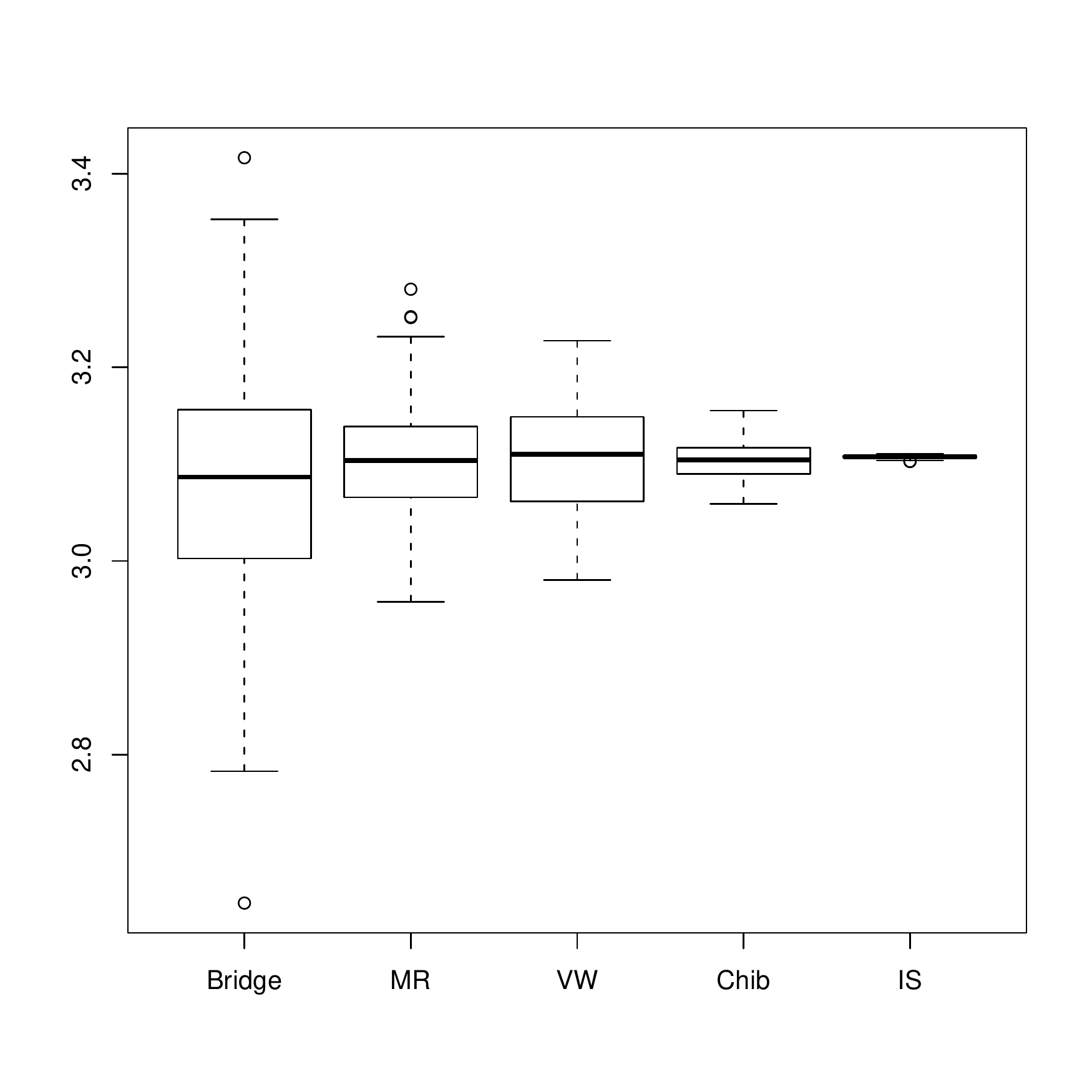}
\caption{\label{fig:bfbsmrvwchiis}
Comparison of the variabilities of five approximations of the Bayes factor evaluating the impact of the diabetes pedigree covariate
upon the occurrence of diabetes in the Pima Indian population, based on a probit modelling. The boxplots are based on $100$ replicas
and the Savage--Dickey representation proposed in the current paper is denoted by MR, while Verdinelli and Wasserman's (1995) version is 
denoted by VW.}
\end{figure}

This probit model is a natural setting for completion by a truncated normal latent variable \citep{albert:chib:1993b}.
We can thus easily implement a Gibbs sampler to produce output from all the posterior distributions considered
in the previous Section.  Besides, in that case, the conditional distribution $\pi_1(\theta|x,\psi,z)$  is
a normal distribution with closed form parameters. It is therefore straightforward to compute the unbiased
estimators \eqref{eq:arrox-mr09} and \eqref{eq:arrox-vw05}. Figure \ref{fig:bfbsmrvwchiis} compares the variation of this approximation
with other standard solutions covered in \cite{marin:robert:2010} for the same example, namely the regular importance
sampling approximation based on the MLE asymptotic distribution, Chib's version based on the same completion,
and a bridge sampling \citep{gelman:meng:1998} solution completing $\pi_0(\cdot)$ with the full
conditional being derived from the conditional MLE asymptotic distribution. The boxplots are all based on 100
replicates of $T=20,000$ simulations. While the estimators \eqref{eq:arrox-mr09} and \eqref{eq:arrox-vw05} are not as
accurate as Chib's version and as the importance sampler in this specific case, their variabilities remain at a reasonable order
and are very comparable. The R code and the reformated datasets used in this Section are available at the following address:
\verb+http://www.math.univ-montp2.fr/~marin/savage/dickey.html+.

\section*{Acknowledgements}
The authors are grateful to H.~Doss and J.~Rousseau for helpful discussions, 
as well as to M.~Kilbinger for bringing the problem to their attention. Comments
from the editorial team were also most useful to improve our exposition of the
Savage--Dickey paradox. The second author also thanks Geoff Nicholls for pointing
out the bridge sampling connection at the CRiSM workshop at the University of
Warwick, May 31, 2010.  This work had been supported by the Agence Nationale de la Recherche (ANR, 212,
rue de Bercy 75012 Paris) through the 2009-2012 project {\sf Big'MC}.


\begin{thebibliography}{18}
\expandafter\ifx\csname natexlab\endcsname\relax\def\natexlab#1{#1}\fi
\expandafter\ifx\csname url\endcsname\relax
  \def\url#1{\texttt{#1}}\fi
\expandafter\ifx\csname urlprefix\endcsname\relax\def\urlprefix{URL }\fi
\providecommand{\eprint}[2][]{\url{#2}}

\bibitem[{Albert and Chib(1993)}]{albert:chib:1993b}
\textsc{Albert, J.} and \textsc{Chib, S.} (1993).
\newblock {B}ayesian analysis of binary and polychotomous response data.
\newblock \textit{J. American Statist. Assoc.}, \textbf{88} 669--679.

\bibitem[{Billingsley(1986)}]{billingsley:1986}
\textsc{Billingsley, P.} (1986).
\newblock \textit{Probability and Measure}.
\newblock 2nd ed. John Wiley, New York.

\bibitem[{Chen et~al.(2000)Chen, Shao and Ibrahim}]{chen:shao:ibrahim:2000}
\textsc{Chen, M.}, \textsc{Shao, Q.} and \textsc{Ibrahim, J.} (2000).
\newblock \textit{{M}onte {C}arlo Methods in {B}ayesian Computation}.
\newblock Springer-Verlag, New York.

\bibitem[{Chib(1995)}]{chib:1995}
\textsc{Chib, S.} (1995).
\newblock Marginal likelihood from the {G}ibbs output.
\newblock \textit{J. American Statist. Assoc.}, \textbf{90} 1313--1321.

\bibitem[{Chopin and Robert(2010)}]{chopin:robert:2010}
\textsc{Chopin, N.} and \textsc{Robert, C.} (2010).
\newblock Properties of evidence.
\newblock \textit{Biometrika}.
\newblock To appear.

\bibitem[{Consonni and Veronese(2008)}]{consonni:veronese:2008}
\textsc{Consonni, G.} and \textsc{Veronese, P.} (2008).
\newblock Compatibility of prior specifications across linear models.
\newblock \textit{Statist. Science}, \textbf{23} 332--353.

\bibitem[{Dickey(1971)}]{dickey:1971}
\textsc{Dickey, J.} (1971).
\newblock The weighted likelihood ratio, linear hypotheses on normal location
  parameters.
\newblock \textit{Ann. Mathemat. Statist.}, \textbf{42} 204--223.

\bibitem[{Gelman and Meng(1998)}]{gelman:meng:1998}
\textsc{Gelman, A.} and \textsc{Meng, X.} (1998).
\newblock Simulating normalizing constants: From importance sampling to bridge
  sampling to path sampling.
\newblock \textit{Statist. Science}, \textbf{13} 163--185.

\bibitem[{Jeffreys(1939)}]{jeffreys:1939}
\textsc{Jeffreys, H.} (1939).
\newblock \textit{Theory of Probability}.
\newblock 1st ed. The Clarendon Press, Oxford.

\bibitem[{Marin and Robert(2010)}]{marin:robert:2010}
\textsc{Marin, J.} and \textsc{Robert, C.} (2010).
\newblock Importance sampling methods for {B}ayesian discrimination between
  embedded models.
\newblock In \textit{Frontiers of Statistical Decision Making and {B}ayesian
  Analysis} (M.-H. Chen, D.~Dey, P.~M{\"u}ller, D.~Sun and K.~Ye, eds.).
  Springer-Verlag, New York.
\newblock To appear, see arXiv:0910.2325.

\bibitem[{Marin and Robert(2007)}]{marin:robert:2007}
\textsc{Marin, J.-M.} and \textsc{Robert, C.} (2007).
\newblock \textit{Bayesian Core}.
\newblock Springer-Verlag, New York.

\bibitem[{{O'H}agan and Forster(2004)}]{ohagan:forster:2004}
\textsc{{O'H}agan, A.} and \textsc{Forster, J.} (2004).
\newblock \textit{Kendall's advanced theory of Statistics: {B}ayesian
  inference}.
\newblock Arnold, London.

\bibitem[{{R Development Core Team}(2008)}]{rmanual}
\textsc{{R Development Core Team}} (2008).
\newblock \textit{R: A Language and Environment for Statistical Computing}.
\newblock R Foundation for Statistical Computing, Vienna, Austria.
\newblock {ISBN} 3-900051-07-0, \urlprefix\url{http://www.R-project.org}.

\bibitem[{Robert(2001)}]{robert:2001}
\textsc{Robert, C.} (2001).
\newblock \textit{The {B}ayesian Choice}.
\newblock 2nd ed. Springer-Verlag, New York.

\bibitem[{Torrie and Valleau(1977)}]{torrie:valleau:1977}
\textsc{Torrie, G.} and \textsc{Valleau, J.} (1977).
\newblock Nonphysical sampling distributions in {M}onte {C}arlo free-energy
  estimation: Umbrella sampling.
\newblock \textit{J. Comp. Phys.}, \textbf{23} 187--199.

\bibitem[{Verdinelli and Wasserman(1995)}]{verdinelli:wasserman:1995}
\textsc{Verdinelli, I.} and \textsc{Wasserman, L.} (1995).
\newblock Computing {B}ayes factors using a generalization of the
  {S}avage--{D}ickey density ratio.
\newblock \textit{J. American Statist. Assoc.}, \textbf{90} 614--618.

\bibitem[{Wetzels et~al.(2010)Wetzels, Grasman and
  Wagenmakers}]{wetzels:grasman:wagenmakers:2010}
\textsc{Wetzels, R.}, \textsc{Grasman, R.} and \textsc{Wagenmakers, E.-J.}
  (2010).
\newblock An encompassing prior generalization of the {S}avage-{D}ickey density
  ratio.
\newblock \textit{Comput. Statist. Data Anal.}, \textbf{54} 2094--2102.

\bibitem[{Zellner(1986)}]{zellner:1986}
\textsc{Zellner, A.} (1986).
\newblock {O}n assessing prior distributions and {B}ayesian regression analysis
  with $g$-prior distribution regression using {B}ayesian variable selection.
\newblock In \textit{{B}ayesian inference and decision techniques: {E}ssays in
  Honor of {B}runo de {F}inetti}. North-Holland / Elsevier, 233--243.

\end{thebibliography}

\end{document}